\documentclass[a4paper,twoside,12pt]{article}
\usepackage[english]{babel}
\usepackage{graphicx}
\usepackage{bbding}
\usepackage[latin1]{inputenc}
\usepackage{kpfonts}
\usepackage{fancyhdr}
\usepackage{amsmath,amssymb,color,epsfig}
\usepackage{mathrsfs}
\usepackage{eufrak}
\usepackage{dsfont}
\usepackage{upgreek}
\usepackage{fancyhdr}
\usepackage[dvipsnames]{xcolor}
\newtheorem{theorem}{Theorem}[section]
\newtheorem{proposition}[theorem]{Proposition}

\newtheorem{remark}[theorem]{Remark}
\newtheorem{lemma}[theorem]{Lemma}

\def\tr{{\rm tr}}

\def\1{\mathds{1}}
\def\e{\varepsilon}

\title{Constrained Eigenvalues Density of Invariant Random Matrices Ensembles }
\date{}
\definecolor{mon}{rgb}{10,126,140}
\begin{document}
\author{Mohamed BOUALI}
\maketitle
\begin{abstract}
We compute exact asymptotic of the statistical density  of random matrices belonging to invariant random matrices ensemble (RMT) orthogonal, unitary
and symplectic ensembles, where all its eigenvalues lie  within the  interval $[\sigma, +\infty[$ or $]-\infty,\tau]$ or $[\sigma,\tau]$. It is found that the density of eigenvalues generically exhibits
an inverse square-root singularity at the location of the barriers. These results generalized the case of Gaussian random matrices ensemble studied in \cite{D}, \cite{S} and \cite{boo}.\\
{\bf Math Subject classification:} 15B52, 15B57, 60B10.
\end{abstract}

\section{Introduction}
Random matrix theory has been successfully applied in various branches of physics and mathematics,
including in subjects ranging from nuclear physics, quantum chaos, disordered systems, and
number theory. Of particular importance are invariant random matrices with density $e^{-\tr(Q(x))}$ with respect the Lebesgue measure where $Q$ is some regular potential. There are three interesting classes of matrices distributed with such density: $(n \times n)$ real symmetric (Gaussian Orthogonal Ensemble (GOE)), $(n \times n)$ complex Hermitian (Gaussian
Unitary Ensemble (GUE)) and $(2n \times 2n)$ self-dual Hermitian matrices ( Gaussian Symplectic Ensemble (GSE)). In
these models the probability distribution for a matrix $X$ in the ensemble is given by
$$p_n(X)\propto \exp\Big(-n\frac\beta 2\tr Q(X)\Big).$$
where $\tr$ is the trace function and $\beta$ is the Dyson index and $n$ is a rescaling factor.

A crucial result in the theory of random matrices is the celebrated Wigner semi-circle law. It states that for regular potential $Q$ and for large
$n$ and on an average, the $n$ eigenvalues of a matrix distributed as the density $p_n(X)$, lie within a finite interval $[-\sqrt 2, \sqrt 2]$, often referred to as the Wigner 'sea'. Within this sea, the statistical density of eigenvalues has a semi-circular form that vanishes at the two edges $-\sqrt 2$, $\sqrt 2$.
$$\rho(\lambda)=\frac1\pi\sqrt{2-\lambda^2}.$$
See for example \cite{LE}, \cite{dei}.
The above result means that, if one looks at the statistical density of eigenvalues of a typical system described by one of the three
ensembles above, for a large enough $n$, it will resemble closely to the Wigner semi-circle law.
From the semi-circle law, we know that on an average half the eigenvalues are positive and half of them are negative.

An important question will be studied in this work, namely  what is the statistical density of eigenvalues which lie  within the  interval $[\sigma, +\infty[$ or $]-\infty,\tau]$ or $[\sigma,\tau]$.

   In this paper we calculate the asymptotic density of eigenvalues in this conditioned ensemble, we will see that it is quite different from the Wigner semi-circle law. We prove the following result.
   The density of eigenvalues in the interval $[\sigma,\tau]$ is given in this way
   $$\int_{\Bbb R}f(x)\mu^{\sigma,\tau}(dx)=\frac1{\pi}\int_{\sigma}^{\tau} f(x)p_{\tau,\sigma}(x)\frac1{\sqrt{(\tau-x)(x-\sigma)}}dx,$$
where $p_{\tau,\sigma}(x)$ is some positive function which will be given later.

Also We prove that as $\tau\rightarrow +\infty$ the density present an inverse square root at $\sigma$ and $p_{\tau,\sigma}(x)\rightarrow p_{\sigma}(x)$, $\mu^{\tau,\sigma}\rightarrow\mu^\sigma$ with
$$\int_{\Bbb R}f(x)\mu^{\sigma}(dx)=\frac1{\pi}\int_{\sigma}^{L(\sigma)} f(x)p_{\sigma}(x)\sqrt{\frac{L(\sigma)-x}{x-\sigma}}dx,$$
where $L(\sigma)$ is some constant which depend on $\sigma$.

 The same hold if one keep $\tau$ fixed and let $\sigma\rightarrow -\infty$ with $\sigma$ replaces by $\tau$.
    Moreover if $\sigma\to-\infty$ and $\tau\to +\infty$ the case of unconstraint eigenvalues, one recovers the Wigner semi-circular law.

   The paper is organized as follows. In the second section we begin by recalling some result about potential theory and equilibrium measure and we enunciate our fundamental equilibrium measure which is the key of the work (theorem \ref{t1}). In section 3 we prove the existence of such measure and we give its density explicitly. Section 4 is dedicated to prove that the measure of theorem \ref{t1} is the limit of the global density of eigenvalues such that the intervals $]-\infty,\sigma]$ or $[\tau,+\infty[$ or $\Bbb R\setminus[\sigma,\tau]$ are devoid of eigenvalues. At the end we argue our study with some examples and numerical simulations.

\section{Preliminary}
Let $\Sigma$ be a closed interval, and $Q$ be a lower semi-continuous function on $\Sigma$. If $\Sigma$ is unbounded we assume that
$$\lim_{|x|\to +\infty}\Big(Q(x)-\log(1+x^2)\Big)=\infty.$$
For given $Q$ and $\Sigma$, we wish to compute the equilibrium measure. We start by some general results.

For any probability measure $\mu$ on $\Sigma$, we define the potential of $\mu$ by: for $x\in\Sigma\setminus{\rm supp}(\mu)$
$$U^\mu(x)=\int_\Sigma\log\frac{1}{|x-y|}\mu(dy),$$
and the energy by
$$E_{Q, \Sigma}(\mu)=\int_\Sigma U^\mu(x)\mu(dx)+\int_\Sigma Q(x)\mu(dx).$$
From the inequality $$|x-y|\leq\sqrt{1+x^2}\sqrt{1+y^2},$$
it can be seen that $E_{Q, \Sigma}(\mu)$ is bounded below. Let
$$E^*_{Q, \Sigma}=\inf_{\mu\in\mathfrak{M}(\Sigma)}E_{Q, \Sigma}(\mu),$$
where $\mathfrak{M}(\Sigma)$ is the set of probability measures on the closed set $\Sigma$.

If $\mu(dx)=f(x)dx$, where $f$ is continuous function with compactly support $\subset\Sigma$, the potential is a continuous function, and $E_{Q, \Sigma}(\mu)<\infty$.
\begin{proposition}---\label{p2}
There is a unique measure $\mu^*\in\mathfrak{M}(\Sigma)$ such that $$E^*_{Q, \Sigma}=E_{Q, \Sigma}(\mu^*),$$
moreover the support of $\mu^*$ is compact.
\end{proposition}
This measure $\mu^*$ is called the {\it equilibrium measure.}

See theorem II.2.3 \cite{F}.
The next proposition is a method to find the equilibrium measure
\begin{proposition}---\label{p1}
Let $\mu\in\mathfrak{M}(\Sigma)$ with compact support.
Assume the potential $U^\mu$ of $\mu$ is continuous and, there is a
constant $C$ such that\\
{\rm (i)} $U^\mu(x)+\frac12Q(x)\geq C$ on $\Sigma$.\\
{\rm (ii)} $U^\mu(x)+\frac12Q(x)= C$ on {\rm supp}$(\mu)$. Then $\mu$ is the equilibrium measure $\mu=\mu^*$.

The constant C is called the (modified) Robin constant. Observe that $$E^*_{Q, \Sigma}=C+\frac12\int_\Sigma Q(x)\mu^*(dx).$$

\end{proposition}
\begin{remark}---
Let $\Sigma'\subset\Sigma$ be a closed interval of $\Sigma$, if we consider the restriction of the function $Q$ initially defined on $\Sigma$ to the closed interval $\Sigma'$ and if the equilibrium measure $\mu$ associate to $(\Sigma, Q)$ satisfies supp$(\mu)\subset\Sigma'$. Then the equilibrium measure for the couple $(\Sigma', Q)$ is $\mu$.
\end{remark}

We saw that if $\Sigma=\Bbb R$, and $Q$ a polynomials of degree $2m$ then the equilibrium measure $\mu_0$ for $(Q,\Sigma)$ is given by Pastur formula
 $$\int_{\Bbb R}f(x)\mu_0(dx)=\frac1{\pi}\int_{a_0}^{b_0} f(x)q(x)\sqrt{(b_0-x)(x-a_0)}dx,$$
 where $q$ is a polynomial of degree $2m-2$ given by
 $$q(x)=\frac1{2\pi}\int_{a_0}^{b_0}\frac{Q'(x)-Q'(t)}{x-t}\frac{dt}{\sqrt{(b_0-t)(t-a_0)}},$$
 The numbers $a_0$ and $b_0$ are given by
 $$\int_{a_0}^{b_0}\frac{Q'(t)dt}{\sqrt{(b_0-t)(t-a_0)}}=0,\;\;\int_{a_0}^{b_0}\frac{tQ'(t)dt}{\sqrt{(b_0-t)(t-a_0)}}=2\pi.$$

We come to our first result. Let $Q$ be a convex polynomial of degree $2m$ with $m\geq 1$ and $\Sigma=[\sigma,\tau]$ with $-\infty\leq\sigma<\tau\leq+\infty$.

{\rm\bf A.} The equilibrium measure of the pair $(Q,\Sigma)$ depend on the parameters $\sigma$ and $\tau$. We denoted it by
$$\mu(\sigma,\tau,dx).$$
Its support is a compact interval $[a,b]$, $a=a(\sigma,\tau)$ and $b=b(\sigma,\tau)$. Its density $h(x):=h(\sigma,\tau; x)$ has singularity at $a$ and $b$.
The following cases are present
\begin{theorem}\label{t1} \
\begin{description}
\item[(1)] { Soft edges in $a$ and $b$}:
$$h(x)=\frac1\pi q(x)\sqrt{(x-a)(b-x)}\chi(x),$$
with $\chi$ the indicator of the interval $[a,b]$. $q$ is a polynomial of degree $2m-2$, given by the Pastur formula
$$q(x)=\frac1{2\pi}\int_{a}^{b}\frac{Q'(x)-Q'(t)}{x-t}\frac{dt}{\sqrt{(b-t)(t-a)}}.$$
\item[(2)] { Hard edge in $a$ and soft edge in $b$}: In that case $a=\sigma$:
$$h(x)=\frac1\pi p(x)\sqrt{\frac{b-x}{x-a}}\chi(x),$$
$p$ is a polynomial of degree $2m-1$ and is given in this way:
$$p(x)=\frac1{2\pi} \int_{a}^{b}\frac{Q'(x)-Q'(t)}{x-t}\sqrt{\frac{t-a}{b-t}}dt+\frac12Q'(x).$$
\item[(3)] { Soft edge in $a$ and hard edge in $b$}: In that case $b=\tau$:
$$h(x)=\frac1\pi s(x)\sqrt{\frac{x-a}{b-x}}\chi(x),$$
$s$ is a polynomial of degree $2m-1$ and is given in this way:
$$s(x)=\frac1{2\pi} \int_{a}^{b}\frac{Q'(x)-Q'(t)}{x-t}\sqrt{\frac{b-t}{t-a}}dt-\frac12Q'(x).$$
\item[(4)] { Hard edges in $a$ and $b$}: then $a=\sigma$, $b=\tau$:
$$h(x)=\frac1\pi r(x)\frac1{\sqrt{(x-a)(b-x)}}\chi(x),$$
where $r$ is a polynomial with degree $2m$, $$r(x)=\frac1{2\pi}\int_{\sigma}^{\tau}\frac{Q'(x)-Q'(t)}{x-t}\sqrt{(t-\sigma)(\tau-t)}dt+\frac12
  \big(\frac{\sigma+\tau}{2}-x\big)Q'(x)+1.$$
\end{description}
\end{theorem}
{\rm\bf B.} Given $\sigma$ and $\tau$. How to determine which case is present and how to find the parameters $a$ and $b$?
\begin{theorem}\label{t22} \
\begin{description}
\item[(i)] For $\sigma=-\infty$ and $\tau=\infty$, let $a=a(-\infty,\infty)=a_0$ and $b=b(-\infty,\infty)=b_0$, and denote $\mu_0$ the equilibrium measure. One can see that we are in the first case (1). If $\sigma\leq a_0$ and $\tau\geq b_0$, then $a=a_0$ and $b=b_0$ and the equilibrium measure is given by $\mu_0$.
\item[(ii)] Fix $\tau=\infty$ and keep the barrier at $\sigma$ moving
\begin{description}
\item[(ii.a)] If $\sigma\leq a_0$ then $\mu(\sigma,\infty,dx)=\mu_0(dx).$
\item[(ii.b)] If $\sigma> a_0$ then $a(\sigma,\infty)=\sigma$, with hard edge in $a(\sigma,\infty)=\sigma$ and soft edge in $b=b(\sigma,\infty)$.
\end{description}
\item[(iii)] We move the barrier $\tau$ in decreasing order from $\infty$.\\
Assume (ii.a) is true
\begin{description}

\item[(iii.a)] If $\tau\geq b_0$ then $\mu(\sigma,\tau,dx)=\mu_0(dx).$
\item[(iii.b)] If $\tau< b_0$ then $b(\sigma,\tau)=\tau$, with soft edge in $a=a(-\infty,\tau)$ and hard edge in $b(\sigma,\tau)=\tau$.

\end{description}
Assume (ii.b) is true
\begin{description}

\item[(iii.c)] If $\tau\geq b(\sigma,\infty)$ then $\mu(\sigma,\tau,dx)=\mu(\sigma,\infty;dx).$ Hard edge in $a(\sigma,\infty)=\sigma$ and soft edge in $b=b(\sigma,\infty)$.
    \item[(iii.d)] If $\tau< b(\sigma,\infty)$, then $a(\sigma,\tau)=\sigma$ and $b(\sigma,\tau)=\tau$. Hard edges in $\sigma$ and $\tau$.
    \end{description}
\end{description}
\end{theorem}

{\bf Proof of theorem \ref{t1}.} We give just the proof of (2) and (4). The others steps follow in the same way.\\
(2).
Let $\mu^\sigma$ be the equilibrium measure of $\Sigma=[\sigma,+\infty[$, $G_\sigma$ be the Cauchy-Stieljes transform of $\mu^\sigma$, $\displaystyle G_\sigma(z)=\int_a^{b}\frac1{z-t}\mu^\sigma(dt).$

$G_\sigma$ is an holomorphic function on  $\Bbb C\setminus[a,b]$ and
$$[G_\sigma]=-2i\pi\mu^\sigma.$$
Moreover, for $x\in[a,b]$
$${\rm Re}G_\sigma(x)=\frac12Q'(x).$$
Let \begin{equation}\label{ww}\widetilde G_\sigma(z)=G_\sigma(z)\sqrt{\frac{z-a}{z-b}}.\end{equation}
The function $\widetilde G _\sigma$ is holomorphic in  $\Bbb C\setminus[a,b]$ and satisfies
$$[\widetilde G_\sigma]=-i\sqrt{\frac{x-a}{b-x}}Q'(x)\chi(x).$$
Since $\widetilde G_\sigma$ goes to $0$ at infinity, using Liouville theorem one gets
\begin{equation}\label{w}\begin{aligned}
\widetilde G_\sigma(z)&=\frac1{2\pi}\int_a^{b}\frac{1}{z-t}Q'(t)\sqrt{\frac{t-a}{b-t}}dt\\
&=-\frac1{2\pi}\int_a^{b}\frac{Q'(z)-Q'(t)}{z-t}\sqrt{\frac{t-a}{b-t}}dt+Q'(z)\frac1{2\pi}\int_a^{b}
\frac{1}{z-t}\sqrt{\frac{t-a}{b-t}}dt.
\end{aligned}\end{equation}
Using equations (\ref{w}) and (\ref{ww}) and lemma \ref{li} in appendix, it follows
$$[G_\sigma]=-2i\sqrt{\frac{b(\sigma)-x}{x-\sigma}} p(x)\chi(x).$$
Hence one deduces the formula as announced, since $[G_\sigma]=-2i\pi\mu^\sigma$.

(4).
Let $\mu^{\sigma,\tau}$ be the equilibrium measure of $\Sigma=[\sigma,\tau]$, and $\displaystyle G_{\sigma,\tau}=\int_a^b\frac1{z-t}\mu^{\sigma,\tau}(dt)$. As in the previous
$$[G_\sigma]=-2i\pi\mu^{\sigma,\tau}.$$
Moreover, for $x\in[a,b]$
$${\rm Re}G_\sigma(x)=\frac12Q'(x).$$
Let \begin{equation}\label{ww}\widetilde G_\sigma(z)=G_\sigma(z)\sqrt{(z-a)(z-b)}.\end{equation}
The function $\widetilde G _\sigma$ is holomorphic in  $\Bbb C\setminus[a,b]$ and satisfies
$$[\widetilde G_\sigma]=i\sqrt{(x-a)(b-x)}Q'(x)\chi(x).$$
Since $\widetilde G_\sigma$ goes to $1$ at infinity, using Liouville theorem one gets
\begin{equation}\label{w}\begin{aligned}
\widetilde G_\sigma(z)&=-\frac1{2\pi}\int_a^{b}\frac{1}{z-t}Q'(t)\sqrt{(t-a)(b-t)}dt+1\\
&=\frac1{2\pi}\int_a^{b}\frac{Q'(z)-Q'(t)}{z-t}\sqrt{(t-a)(b-t)}dt-Q'(z)\frac1{2\pi}\int_a^{b}
\frac{1}{z-t}\sqrt{(t-a)(b-t)}dt+1.
\end{aligned}\end{equation}
using equations (\ref{w}) and (\ref{ww}) and lemma \ref{li} in appendix, it follows
$$[G_\sigma]=-2i\frac1{\sqrt{(b-x)(x-a)}} r(x)\chi(x).$$
Hence one deduces the formula as announced.

{\bf Proof of theorem \ref{t22}.}\\
{\it Proof of (ii).}\\
(ii.b).
Let $\Sigma_\sigma=[\sigma,+\infty[$ with $\sigma\geq a_0$ and $Q$ any function satisfying the hypotheses as above. Let $\mu^\sigma(dx):=\mu(\sigma,\infty,dx)$ be the measure defined as in (2) of theorem \ref{t1}, with $b=b(\sigma,\infty)=b(\sigma)$ given in the proposition \ref{pp1} in appendix,
$$p(x)=\frac1{2\pi} \int_{\sigma}^{b(\sigma)}\frac{Q'(x)-Q'(t)}{x-t}\sqrt{\frac{t-\sigma}{b(\sigma)-t}}dt+\frac12Q'(x).$$
Remark that from lemma \ref{li} and proposition \ref{pp1} in appendix, one gets $$\displaystyle p(\sigma)=\int_\sigma^{b(\sigma)}\frac{Q'(t)}{\sqrt{(t-\sigma)(b(\sigma)-t)}}=\varphi(\sigma,b(\sigma))\geq 0.$$ Moreover a simple computation give
 $$p(x)=p(\sigma)+(x-\sigma)q(x),$$
 where $$q(x)=\frac1{2\pi}\int_{\sigma}^{b(\sigma)}\frac{Q'(x)-Q'(t)}{x-t}\frac{dt}{\sqrt{(b(\sigma)-t)(t-\sigma)}}.$$
 By convexity of the polynomial $Q$ and the fact that $p(\sigma)\geq 0$, one gets for all $x\geq\sigma$, $p(x)\geq 0$.

Expanding the function $\widetilde G_\sigma$ near infinity, it follows
$$\widetilde G_\sigma(z)=\sum_{n=0}^{+\infty}\frac{a_n}{z^{n+1}},$$
where $\displaystyle a_n=\frac1{2\pi}\int_{\sigma}^{b(\sigma)}t^nQ'(t)\sqrt{\frac{t-\sigma}{b(\sigma)-t}}dt,$
Since $$\frac1{2\pi}\int_{\sigma}^{b(\sigma)}Q'(t)\sqrt{\frac{t-\sigma}{b(\sigma)-t}}dt=-\frac1{2\pi}\psi(\sigma,b(\sigma))+1=1,$$
which follow from proposition \ref{pp1} in appendix.
Thus $$\widetilde G_\sigma(z)=\frac1z+o(\frac1z),$$
moreover $$G(z)=\sqrt{\frac{z-\sigma}{z-b(\sigma)}}\widetilde G_\sigma(z),$$
hence,
$$G(z)=\frac1z+o(\frac1z).$$

Which prove that $\mu(\sigma,\infty;dx)$ is a probability.

 It remains to prove that $\mu(\sigma,\infty;dx)$ is the equilibrium measure.
Since for $a=\sigma$ and $b=b(\sigma)$,
$$G_\sigma(z)=-\sqrt{\frac{z-b(\sigma)}{z-\sigma}}p(z)+\frac12Q'(z).$$
Moreover $Q$ is holomorphic, hence $p$ is holomorphic too. It follows that

$$\lim_{y\to 0,\,y>0}{\rm Re}G_\sigma(x+iy)=\left\{\begin{aligned}&-\sqrt{\frac{b(\sigma)-x}{\sigma-x}}p(x)+\frac12Q'(x)\;\; if\;\;x\geq b(\sigma),\\
&\quad\frac12Q'(x)\;\qquad\qquad\qquad\quad{if}\,\sigma\leq x\leq b(\sigma).
\end{aligned}\right.$$
An easy computation shows that $$\frac{d}{dx}U^{\mu^\sigma}(x)=-\lim_{y\to 0,\,y>0}{\rm Re}G_\sigma(x+iy).$$
Hence

$$2\frac{d}{dx}U^{\mu^\sigma}(x)+Q'(x)=\left\{\begin{aligned}&2\sqrt{\frac{b(\sigma)-x}{\sigma-x}}p(x)\quad if\;\;x\geq b(\sigma),\\
&0\qquad\qquad\qquad\;{if}\,\sigma\leq x\leq b(\sigma).
\end{aligned}\right.$$
We have proved in the beginning  that $p(x)\geq0$ for all $x\geq\sigma\geq a_0$, hence
 $$2U^{\mu^\sigma}(x)+Q(x)\left\{\begin{aligned}&\geq C\qquad\qquad if\;x\geq b(\sigma),\\
&=C\qquad\qquad\;{if}\,\sigma\leq x\leq b(\sigma).
\end{aligned}\right.$$
(ii.a). For $\sigma\leq a_0$, The equilibrium measure is $\mu_0(dx)$. In fact, if there is some $a\in[\sigma,a_0[$ such that $p(x)\geq 0$ for all $x\in[a,a_0[$,
where $$p(x)=\frac1{2\pi} \int_{a}^{b(a)}\frac{Q'(x)-Q'(t)}{x-t}\sqrt{\frac{t-a}{b(a)-t}}dt+\frac12Q'(x).$$
 Then by the equation
$$p(x)=p(a)+(x-a)q(x),$$
one gets for $x\in[a,a_0[$ $$p(a)+(x-a)q(x)\geq 0.$$
As $x\to a$, we obtain $$p(a)=\varphi(a,b(a))\geq 0.$$
This contradicts the minimality of $a_0$ see proposition \ref{pp1} in appendix.
This complete the proof.\\
{\it Prove of  (iii).} We prove just (iii.c) and (iii.d), the rest hold with the same lines as in (ii) with little modification.\\
$1)$ $a_0<\sigma$ and $\Sigma=[\sigma,\tau]$.\\
(iii.c) First case $\tau<b(\sigma)$.
The density of the measure is
$$\mu^{\sigma,\tau}(dx)=h(x,\sigma,\tau)=\frac1\pi\frac{r(x)}{\sqrt{(\tau-x)(x-\sigma)}}\chi_{[\sigma,\tau]}(x)dx,$$
with $$r(x)=\frac1{2\pi}\int_{\sigma}^{\tau}\frac{Q'(x)-Q'(t)}{x-t}\sqrt{(t-\sigma)(\tau-t)}dt-\frac12\Big(x-\frac{\sigma+\tau}{2})\Big)Q'(x)+1.$$
Since by some algebra one gets $$r(x)=r(\tau)+(\tau-x)p(x),$$
where $$p(x)=\frac1{2\pi} \int_{\sigma}^{\tau}\frac{Q'(x)-Q'(t)}{x-t}\sqrt{\frac{t-\sigma}{\tau-t}}dt+\frac12Q'(x),$$
moreover $p(x)\geq 0$ for all $x\geq\sigma$. It remains to prove that $r(\tau)\geq 0$.

A simple computation gives $$r(\tau)=-\frac1{2\pi}\psi(\sigma,\tau),$$
hence $$r(\tau)=-\frac1{2\pi}h(\sigma,\tau)(\tau-\sigma).$$
We saw that the function $x\mapsto h(\sigma,x)$ increases for all $x>\sigma$, and $b(\sigma)\geq\tau>\sigma$, it follows that
$$r(\tau)=-\frac1{2\pi}h(\sigma,\tau)(\tau-\sigma)\geq-\frac1{2\pi}h(\sigma,b(\sigma)(\tau-\sigma)=0.$$

If $$G_{\sigma,\tau}(z)=\int_\sigma^{\tau}\frac1{z-t}\mu^{\sigma,\tau}(dt),$$
then as in the previous, one has
$$G_{\sigma,\tau}(z)=\frac1{\sqrt{(z-\sigma)(z-\tau)}}S(z),$$
where $$S(z)=\frac1{2\pi}\int_{\sigma}^{\tau}\frac{Q'(z)-Q'(t)}{z-t}\sqrt{(t-\sigma)(\tau-t)}dt-\frac12\Big(z-\frac{\sigma+\tau}{2}
-\sqrt{(z-\sigma)(z-\tau)}\Big)Q'(z)+1.$$
Since for $z\in\Bbb C\setminus[\sigma,\tau]$, $$S(z)=\frac1{2\pi}\int_{\sigma}^{\tau}\frac{-Q'(t)}{z-t}\sqrt{(t-\sigma)(\tau-t)}dt+1,$$
thus $$\lim_{|z|\to +\infty}S(z)=1,$$
and $$\lim_{|z|\to+\infty}|z|G(z)=\lim_{|z|\to +\infty}S(z)=1.$$
It follows that $\mu^{\sigma,\tau}$ is a probability measure.
Moreover
$$\lim_{y\to 0,\,y>0}{\rm Re}G_\sigma(x+iy)=\frac12Q'(x),\quad\forall x\in[\sigma,\tau].$$
Hence $$2\frac{d}{dx}U^{\mu^{\sigma,\tau}}(x)+Q'(x)=0\quad\forall\;x\in [\sigma,\tau].$$
This shows that $\mu^{\sigma,\tau}$ is the equilibrium measure.

(iii.d) Second case $\tau> b(\sigma)$.
In this case we saw that $\mu(\sigma,\infty,dx)$ is a probability measure on $\Sigma=[\sigma,+\infty[$ with support $[\sigma,b(\sigma)]$, hence it is a probability on $[\sigma,\tau]$. Moreover for all $\tau\geq b(\sigma)$

$$\lim_{y\to 0,\,y>0}{\rm Re}G_\sigma(x+iy)=\left\{\begin{aligned}&-\sqrt{\frac{x-b(\sigma)}{x-\sigma}}p(x)+\frac12Q'(x)\;\; if\;\tau\geq x\geq b(\sigma),\\
&\quad\frac12Q'(x)\;\qquad\qquad\qquad\quad{if}\,\sigma\leq x\leq b(\sigma).
\end{aligned}\right.$$
Hence $\mu(\sigma,\infty,dx)$ is an equilibrium measure on $[\sigma,\tau]$. By unicity one gets $\mu^{\sigma,\tau}(dx)=\mu(\sigma,\infty,dx)$.

One can prove (iii.d) in this way. Assume there is some $b\in]b(\sigma),\tau]$ such that the measure $\displaystyle\mu(\sigma,b,dx)=\frac1{\pi}\frac{r(x)}{\sqrt{(x-\sigma)(b-x)}}dx$ define a probability measure on $[\sigma,b]$. Then as in the previous, for all $x\in[\sigma,b]$
$$r(x)=r(b)+(b-x)p(x),$$
with $$p(x)=\frac1{2\pi} \int_{\sigma}^{b}\frac{Q'(x)-Q'(t)}{x-t}\sqrt{\frac{t-\sigma}{b-t}}dt+\frac12Q'(x),$$
We saw that $\displaystyle r(b)=-\frac1{2\pi}\psi(\sigma,b)$.
In other hand, the measure $\mu(\sigma,b,dx)$ is a probability, hence
$\displaystyle-\frac1{2\pi}\psi(\sigma,b)\geq 0$, and $\displaystyle h(\sigma,b)=\frac{\psi(\sigma,b)}{b-\sigma}\leq 0$ which gives $h(\sigma,b)\leq 0$. Furthermore, the function $x\mapsto h(\sigma,x)$ increases strictly, hence
$$0=h(\sigma,b(\sigma))<h(\sigma,b)\leq 0.$$
Which gives a contradiction. This prove that the support of the measure is $[\sigma,b(\sigma)]$ and
$$r(x)=-\frac1{2\pi}\psi(\sigma,b(\sigma))+(b(\sigma)-x)p(x)=(b(\sigma)-x)p(x),$$
hence $$\mu(\sigma,\tau,dx)=\mu(\sigma,\infty,dx).$$

\begin{remark} To see how to pass from the hard edge in $\tau$ to the soft edge in $b(\sigma)$, for $\tau< b(\sigma)$, one can write $r(x)$ in the following sense
$$r(x)=r(\tau)+(\tau-x)p(x),$$
since $r(\tau)=-\frac1{2\pi}\psi(\sigma,\tau)$. As $\tau\to b(\sigma)$, one gets $r(b(\sigma))=0$ see proposition \ref{pp1} in appendix, hence
$$r(x)=(b(\sigma)-x)p(x),$$
and the density $$\frac1{\pi}\frac{r(x)}{\sqrt{(\tau-x)(x-\sigma)}},$$
becomes $$\frac1{\pi}p(x)\sqrt{\frac{b(\sigma)-x}{x-\sigma}}.$$
The same hold in all the others cases.
\end{remark}
\begin{remark} One can change the conditions on the potential $Q$, to be a sufficiently regular and convex function with $\displaystyle\lim_{|x|\to +\infty}Q(x)-\log(1+x^2)=+\infty$  and $\displaystyle\lim_{x\to \pm\infty} Q'(x)=\pm\infty$.
\end{remark}
\section{Constrained eigenvalues density}
We consider the invariant random matrices ensemble with Dyson index $\beta = 1,2,4$, corresponding
to real, complex, and quaternion entries, respectively. The
probability distribution of the entries is given by
$$\Bbb{P}_{n}(dX)=\frac1{C_n}\exp\Big(-n\frac{\beta}{2}\tr Q(X)\Big)dX.$$
where $Q$ is a convex polynomial with even degree $2m$, $m\geq 1$, $C_n$ is a normalizing constant and $dX$ is the Lebesgue measure on the space $H_n=Herm(n,\Bbb F)$ of hermitian matrices with respectively real, complex or quaternion coefficients $\Bbb F=\Bbb R$, $\Bbb C$, or $\Bbb H$. Consequently the joint probability density of eigenvalues
is given by
$$\Bbb{P}_{n}(d\lambda_1,...,d\lambda_n)=\frac1{C_n}e^{-n\frac{\beta}{2}
\sum\limits_{i=1}^nQ(\lambda_i)}|\Delta(\lambda)|^{\beta}d\lambda_1...d\lambda_n,$$
where $\Delta(\lambda)=\prod_{i<j}(\lambda_i-\lambda_j)$ is the Vandermonde determinant, and $n$ is a normalizing factor.

More generally one can consider this probability with an arbitrary real number $\beta>0$.

 For  $\sigma,\tau\in[-\infty,+\infty]$, $\sigma<\tau$, consider $$\Omega_{n, \sigma,\tau}=\Big\{X\in H_n\mid\,\lambda_{\rm min}(X)\geq \sigma\; {\rm and}\;\lambda_{\rm \max}(X)\leq\tau\Big\},$$ the subset of Hermitian matrices for which all its eigenvalues are in $[\sigma,\tau]$.

 We wish to study $\Bbb{P}_{n}(\Omega_{n, \sigma,\tau})$, the probability for a hermitian matrix $X\in H_n$, to have all its eigenvalues in $\Omega_{n, \sigma,\tau}$. It is the probability that all the eigenvalues lies in the interval $[\sigma, \tau]$ that is
 $$\Bbb{P}_{n}(\Omega_{n, \sigma,\tau})=\frac1{C_n}\int_{[\sigma,\tau]^n}e^{-n\frac{\beta}{2}
\sum\limits_{i=1}^nQ(\lambda_i)}|\Delta(\lambda)|^{\beta}d\lambda_1...d\lambda_n.$$

Let $\nu^{\sigma,\tau}_{n}$ be the probability measure defined on $[\sigma,\tau]$ by : for all continuous functions $f$
$$\int_{[\sigma,\tau]}f(x)\nu^{\sigma,\tau}_{n}(dx)=\int_{[\sigma,\tau]^n}\frac1n\sum_{i=1}^nf(\lambda_i)\Bbb{P}_{n}(d\lambda_1,d\lambda_2,...,d\lambda_n),$$
which means that $$\nu^{\sigma,\tau}_{n}=\Bbb{E}_{n}\big(\frac1n\sum_{i=1}^n\delta_{\lambda_i}\big),$$
where $\Bbb{E}_{n}$ is the expectation with respect to the measure $\Bbb{P}_{ n}$.

As $n$ goes to infinity we prove, the measure $\nu^{\sigma,\tau}_{n}$ converges to some probability measure, which is the statistical density of the eigenvalues in the interval $[\sigma,\tau]$.
In other word one has the following theorem
\begin{theorem}---\label{t2} There exist a unique $a=a(\sigma,\tau)$, and $b=b(\sigma, \tau)$, such that, the measure $\nu^{\sigma,\tau}_{n}$ converges for the tight topology to the probability measure $\mu(\sigma,\tau,dx)$ of theorem \ref{t1}, with support $[a,b]$, and density
$h(\sigma, \tau,x)$.
This means for all continuous functions $\varphi$ on $[\sigma,\tau]$,
$$\lim_{n\to\infty}\int_\sigma^\tau\varphi(x)\nu_{n}^{\sigma,\tau}(dx)=\int_\sigma^\tau\varphi(x)\mu(\sigma,\tau,dx).$$
\end{theorem}

{\it Examples.} \\
1) The density of eigenvalues such that all eigenvalues lie in the interval $[\sigma,+\infty[$ is given by the measure $\mu^\sigma(dx):=\mu(\sigma,\infty,dx)$, see (2) theorem \ref{t1}.\\
2)  The density of eigenvalues in $\Bbb R$ (unconstrained condition) is given by Pastur formula $\mu(-\infty,\infty,dx)$, see (1) theorem \ref{t1}.

\subsection{Proof of theorem \ref{t2}}
 In the rest of the section we prove theorem \ref{t2}, for this purpose we need some preliminary results.

Let $K_n$ be the function on $\Sigma^n=[\sigma,\tau]^n$, defined by
$$K_n(x)=\sum_{i\neq j}k_n(x_i,x_j),$$
where
$$k_n(x,y)=\log\frac{1}{|x-y|}+\frac12Q(x)+\frac12Q(y).$$
The function $K_n$ is bounded below, moreover if $\sigma=-\infty$ or $\tau=+\infty$, $\lim_{|x|\to+\infty}K_n(x)=+\infty$, it follows that $K_n$ attaint it minimum at a point say, $x^{(n, \sigma,\tau)}=(x_1^{(n, \sigma,\tau)},...,x_n^{(n, \sigma,\tau)} )$.
Let $$\delta^{\sigma,\tau}_n=\frac1{n(n-1)}\inf_{\Sigma^n} K_n(x),$$
and $$\rho^{\sigma,\tau}_n=\frac1n\sum_{i=1}^n\delta_{x^{(n, \sigma,\tau)}_i}.$$

For a probability measure $\mu$ on $\Sigma$, consider the energy
$$E_{\delta, \sigma}(\mu)=\int_{\Sigma^2}\log\frac{1}{|s-t|}\mu(ds)\mu(dt)+\int_{\Sigma}Q(s)\mu(ds).$$
and $$E^*_{\delta, \sigma}=\inf_{\mu}E_{\delta, \sigma}(\mu),$$
where the minimum is taken over all compactly support measures with support in $\Sigma$.
Moreover, defined the scaled density
 $$\Bbb{P}_{ n}(dx)=\frac1{Z_n}e^{-n\frac{\beta}{2}
\sum\limits_{i=1}^nQ(x_i)}|\Delta(x)|^{\beta}dx_1...dx_n.$$
where $Z_n$ is a normalizing constant.
\begin{proposition}---\label{p3} \begin{description}
\item[(1)] $\displaystyle\lim_{n\to\infty}\delta^{\sigma,\tau}_n=E^*_{\sigma,\tau}=E_{\sigma}(\mu^{\sigma,\tau}).$
\item[(2)] The measure $\rho^\sigma_n$ converge for the tight topology to the equilibrium measure $\mu^{\sigma,\tau}$.
\item[(3)] $\displaystyle\lim_{n\to\infty}-\frac1{n^2}\log Z_n=\frac\beta2E^*_{ \sigma,\tau}$.
\end{description}
\end{proposition}
{\bf Proof.}---\\
{\bf Step 1 and 2:}
For a probability measure $\mu$,
$$\int_{\Bbb R^n_+}K^\sigma_n(x)\mu(dx_1)...\mu(dx_n)=n(n-1)\int_{\Sigma^2}\log\frac{1}{|x-y|}\mu(dx)\mu(dx)+n(n-1)\int_{\Sigma}Q(x)\mu(dx),$$
hence $$\delta_n^{\sigma,\tau}\leq E_{\sigma,\tau}(\mu).$$
For $\mu=\nu^{\sigma,\tau}$,
 \begin{equation}\label{34}\delta_n^{\sigma,\tau}\leq E^*_{\sigma,\tau}.\end{equation}
Moreover $$K_n(x^{(n, \sigma,\tau)})=\sum_{i\neq j}k_n\Big(x_i^{(n, \sigma,\tau)},x_j^{(n, \sigma,\tau)}\Big)\geq(n-1)\frac12\Big(\sum_{i=1}^nh\Big(x_i^{(n, \sigma,\tau)}\Big)+
\sum_{i=1}^nh\Big(x_i^{(n,\sigma,\tau)}\Big)\Big),$$
where $\displaystyle h(x)=Q(x)-\log(1+x^2),$
Since $$\int_{\Sigma}h(t)\rho^{\sigma,\tau}_n(dt)=\frac1n\sum_{i=1}^nh\Big(x_i^{(n, \sigma,\tau)}\Big),$$
it follows that,
$$\int_{\Sigma}h(t)\rho^{\sigma,\tau}_n(dt)\leq\delta^{\sigma,\tau}_n\leq  E^*_{\sigma,\tau}.$$
Using the fact that $\lim_{x\to \infty}h(x)=+\infty,$ then by the Prokhorov criterium there is some subsequence $\rho_{n_k}$, which convergent to $\rho$ for the tight topology.

For $\ell\geq 0$, let $k_n^{ \ell}(x,y)=\inf(k_n(x,y),\ell)$,
defined $$E^\ell_{\sigma,\tau}(\mu)=\int_{\Sigma^2}k^{ \ell}_n(x,y)\mu(dx)\mu(dy),$$

\begin{equation}\label{35}E^\ell_{ \sigma,\tau}(\rho^{\sigma,\tau}_{n_k})\leq\delta_n^{\sigma,\tau}+\frac{\ell}{n}\leq E^*_{\sigma,\tau}+\frac{\ell}{n}.\end{equation}
The cut kernel $k^{\sigma, \ell}_{\alpha\pm\e}(x,y)$ is bounded and continuous, and the probability measure
 $\rho^{\sigma,\tau}_{n_k}$ converge tightly to $\rho^{\sigma,\tau}$, hence
  $\displaystyle\lim_{k\to+\infty}E_{\sigma,\tau}(\rho^{\sigma,\tau}_{n_k})=E_{\sigma,\tau}(\rho^{\sigma,\tau})$.
As $\ell$ goes to $+\infty$, by the monotone convergence theorem one obtains
$$E_{ \sigma,\tau}(\rho^{\sigma,\tau})\leq E^*_{\sigma,\tau}.$$

By the definition of the equilibrium measure we obtains $E^*_{\sigma,\tau}=E_{\sigma,\tau}(\rho^{\sigma,\tau})=E_{\sigma,\tau}(\mu^{\sigma,\tau})$, it follows by unicity of the equilibrium measure that $\rho^{\sigma,\tau}=\mu^{\sigma,\tau}$. Which means the only possible limit for a subsequence of $\rho_n^{\sigma,\tau}$ is $\mu^{\sigma,\tau}$, hence the sequence $\rho_n^{\sigma,\tau}$ it self converge to $\mu^{\sigma,\tau}$.
Moreover from equation (\ref{35}) one gets
 $$\lim\limits_n\delta^{\sigma,\tau}_n=E^*_{\sigma,\tau}.$$
{\bf Step 3:}
We saw for every $x\geq 0$, $K^\sigma_n(x)\geq n(n-1)\delta^{\sigma,\tau}_n,$ hence
$${Z}_n\leq e^{-\frac\beta 2n(n-1)\delta^{\sigma,\tau}_n}c^n ,$$
where $\displaystyle c=\int_0^{+\infty}e^{-\frac\beta 2Q(x)}dx$,
hence $$\frac1{n^2}\log{Z}_n\leq -\frac\beta2\frac{n-1}{n}\delta^{\sigma,\tau}_n+\frac1n\log(c).$$
Then
$$\limsup_n\frac1{n^2}\log{Z}_n\leq -\frac\beta 2E^*_{\sigma, \tau}.$$
Furthermore $${Z}_n\geq \int_{\Bbb R^n}e^{-\frac\beta 2K_n(x)-\frac\beta 2Q(x)-\sum\limits_{i=1}^n\log h(\sigma,\tau,x_i)}\prod_{i=1}^n\mu^{\sigma,\tau}(dx_i),$$
Applying Jensen's inequality we obtain

$${Z}_n\geq\exp\int_{\Bbb R^n}\Big(-\frac\beta 2K_n(x)-\frac\beta 2Q(x)-\sum\limits_{i=1}^n\log h(\sigma,\tau,x_i)\Big)\prod_{i=1}^n\mu^{\sigma,\tau}(dx_i),$$
hence
$${Z}_n\geq e^{-\frac\beta 2\big(n(n-1)E^*_{\sigma,\tau}\big)}\exp\Big(-\frac\beta 2n\int_a^bQ(x)h(\sigma,\tau, x)dx\Big)\exp\Big(-n\int_a^bh(\sigma,\tau, x)\log h(\sigma,\tau, x)\Big)dx.$$
The function $x\mapsto h(\sigma,\tau, x)\log h(\sigma, \tau, x)$ is continuous on $[a,b]$ $(a:=a(\sigma,\tau); b:=b(\sigma,\tau))$. Hence
$$\liminf_n\frac1{n^2}\log{Z}_n\geq -\frac\beta 2E^*_{\sigma,\tau},$$
and the conclusion hold
 $$-\frac\beta 2E^*_{ \sigma,\tau}\leq\liminf_n\frac1{n^2}\log{Z}_n\leq\limsup_n\frac1{n^2}\log{Z}_n\leq -\frac\beta 2E^*_{\sigma,\tau}.$$
 {\bf Proof of theorem \ref{t2}}---The proof of the theorem follows the proof in (\cite{F}, theorem IV.5.1).
\section{Examples.}
{\it Example 1.}

If one considers the Gaussian invariant ensemble, for such ensemble $Q(x)=x^2$.
The solutions of the two equations $\varphi(b_0,a_0)=\psi(b_0,a_0)=0$ are $(a(-\infty,\infty),b(-\infty,\infty))$ see proposition \ref{pp1} in appendix, one obtains $a(-\infty,\infty)=a_0=-\sqrt 2$, and  $b(-\infty,\infty)=b_0=\sqrt 2$, moreover
 $$r(x)=\frac1{\pi}\int_{\sigma}^{\tau}\sqrt{(t-\sigma)(\tau-t)}dt+x(-x+\frac{\sigma+\tau}{2})+1.$$
{\it Two Hard edges.} By a simple computation  it yields the density of the measure  $\mu^{\sigma,\tau}$: for $-\sqrt 2\leq\sigma<\tau\leq\sqrt 2$,
\begin{equation}\label{eq001} h(\sigma,\tau,x)=\frac1{\pi}\frac1{\sqrt{(\tau-x)(x-\sigma)}}\Big(\frac{(\sigma-\tau)^2}{8}+1+x\frac{\sigma+\tau}{2}-x^2\Big),\end{equation}

 {\it Two Soft edges.} For $\sigma\leq -\sqrt 2$ and $\tau\geq\sqrt 2$ one gets two soft edges in $-\sqrt 2,\sqrt 2$ and  $\mu^{\sigma ,\tau}=\mu^{-\sqrt 2 ,\sqrt 2}$ is the semicircle law with density
 $$h(\sigma,\tau,x)=h(-\infty,\infty,x)=\frac1\pi\sqrt{2-x^2},$$
 $h(-\infty,\infty,x)$ is the density of unconstrained eigenvalues.

 {\it Hard edge in $\sigma$ and Soft edge in $b(\sigma)$.} For $-\sqrt 2<\sigma<\sqrt 2$ and $\tau>\sqrt 2$, one gets a hard edge in $\sigma$ and a soft edge in $b(\sigma)<\tau$ and from proposition \ref{pp1} appendix we obtain $\displaystyle b(\sigma)=\frac23\Big(\frac\sigma2+\sqrt{\sigma^2+6}\Big)$ and the density with support $[\sigma,b(\sigma)]$ is given in this way
 $$h(\sigma,\tau,x)=h(\sigma,\infty,x)=\frac1{2\pi}\sqrt{\frac{b(\sigma)-x}{x-\sigma}}\Big(2x+b(\sigma)-\sigma\Big)\chi(x),$$
The density $h(\sigma,\infty,x)$ represents the density of eigenvalues of Gaussian hermitian random matrices to have all its eigenvalues in $[\sigma,+\infty[$. This cases agree with results of Dean-Majumdar see for instance \cite{D} and \cite{S}.

{\it Soft edge in $a(\tau)$ and Hard edge in $\tau$.} For $\sigma<-\sqrt 2$ and $\tau<\sqrt 2$, we obtain a hard edge in $\tau$ and a soft edge in $\displaystyle a(\tau)=\frac23\big(\frac \tau2-\sqrt{\tau^2+6}\big)$ and the density with support $[a(\tau),\tau]$ is
$$h(\sigma,\tau,x)=h(-\infty,\tau,x)=\frac1{2\pi}\sqrt{\frac{x-a(\tau)}{\tau-x}}\Big(a(\tau)-\tau-2x\Big)\chi(x),$$
$h(-\infty,\tau,x)$ is the density such that all eigenvalues of Gaussian hermitian matrices lie in $]-\infty,\tau]$.

 {\it Example 2.} Let $Q(x)=x^2-\mu_n\log|x|$ with $\mu_n$ a nonnegative sequence of real numbers. It has been proved in \cite{boo} that, if $\displaystyle\lim_{n\to\infty}\frac{\mu_n}n=\alpha$, the density of eigenvalues to be all in the half line $[\sigma,+\infty[$ is
 $$f_{\alpha, a}(x)=\frac{1}{2\pi}\sqrt{\frac{b-x}{x-a}}\Big(2x+b-a-2\alpha\sqrt{\frac{a}{b}}\frac1x\Big)\chi(x),$$
   where $a\geq\sigma$ and $b>a$ are the unique solutions of the following equations
   \begin{equation}\label{01}b+a-\frac{2\alpha}{\sqrt {ab}}=0,\;\;\frac{3}{4}(b-a)^2+a(b-a)+2\alpha\frac{\sqrt a}{\sqrt b}-2\alpha-2=0,\end{equation}
   and $\chi$ is the characteristic function of the interval $[a,b]$. One can see that the two previous equations are just the following $\varphi(a,b)=0$ and $\psi(a,b)=0$ as in proposition\ref{pp1} in appendix.

 If $\alpha>0$, for all $\sigma>0$ the previous equations admit a unique solutions $a(\sigma,\alpha)\geq\sigma$ and $b(\sigma,\alpha)>a(\sigma,\alpha)$. It has been proved that there exist some critical value $a_c(\alpha)$ and $b_c(\alpha)$ as in theorem \ref{t22} $(a_0,b_0)$, which is a transition point from a hard edge in $a$ to a soft edge in $a_c$.

 For $\alpha>0$ then $a_c>0$ and for $ 0\leq\sigma\leq a_c$, the density of eigenvalues to be all in the interval $[\sigma,+\infty[$ is
 \begin{equation}\label{e1}f_\alpha(x)=\sqrt{(b_c-x)(x-a_c)}\Big(1+\frac\alpha{\sqrt {a_c b_c}}\frac1x\Big)\chi(x).\end{equation}

  As $\alpha\to 0^+$ we find from equation (\ref{01}) the case of Gaussian random unitary ensemble. As $\alpha\to 0^+$ then $a_c(\alpha)\to-\sqrt 2$ and $b_c(\alpha)\to\sqrt 2$ and one recovers the Wigner semicircle law.

  Also in that case one can find the density of eigenvalues to lie within in the interval $[0,+\infty[$ which is given by
  $$f(x)=\frac{1}{2\pi}\sqrt{\frac{b-x}{x}}\Big(2x+b\Big)\chi(x),$$
  with $b=\frac23\sqrt 6$ and $\chi$ the indicator function of the interval $[0,\frac23\sqrt 6]$.

{\bf Approximation Density.} Here we give a numerical simulation of the density of eigenvalues for the case where $Q(x)=x^2-\mu_n\log|x|$.

Recall that for $\lim_{n\to+\infty}\frac{\mu_n}{n}=\alpha$, the density of positive eigenvalues is $\displaystyle f_{\alpha}(x)$ of equation (\ref{e1}).

Let $$f_n(x)=\frac1{\sqrt n}\sum_{k=0}^{n-1}\varphi^{\mu_n}_k(\sqrt nx)^2,$$
where $\displaystyle\varphi^{\mu_n }_k(x)=\frac1{\sqrt {d_{k,n}}}H^{\mu_n}_k(x)x^{\mu_n}e^{-\frac{x^2}{2}},$  and $H_k^{\mu_n}$ is the truncated orthogonal Hermite polynomial on the positive real axis, which satisfies $$\int_0^{+\infty}H_k^{\mu_n}(x)H_m^{\mu_n}(x)x^{2\mu_n}e^{-x^2}dx=0,\quad {\rm for}\; m \neq k,$$
and $$\int_0^{+\infty}(H_k^{\mu_n}(x))^2x^{2\mu_n}e^{-x^2}dx=d_{k, n}.$$
It has been proved in theorem 5.1 in \cite{boo}, that as n go to $+\infty$, the density $f_n$ converge tightly to the density $f_\alpha$ where $\displaystyle\alpha=\lim_{n\to\infty}\frac{\mu_n}n$.

{\bf First case:} $n=7, \mu_7=0$, hence $\alpha=0$.
$$\displaystyle f_0(x)=\frac1{2\pi}\sqrt{\frac{\frac23\sqrt 6-x}{x}}\big(2x+\frac23\sqrt 6\big), {\rm with\;\; support}\;\; [0,\frac23\sqrt 6].$$

\begin{figure}[h]
\centering\scalebox{0.5}{\includegraphics[width=17cm, height=12cm]{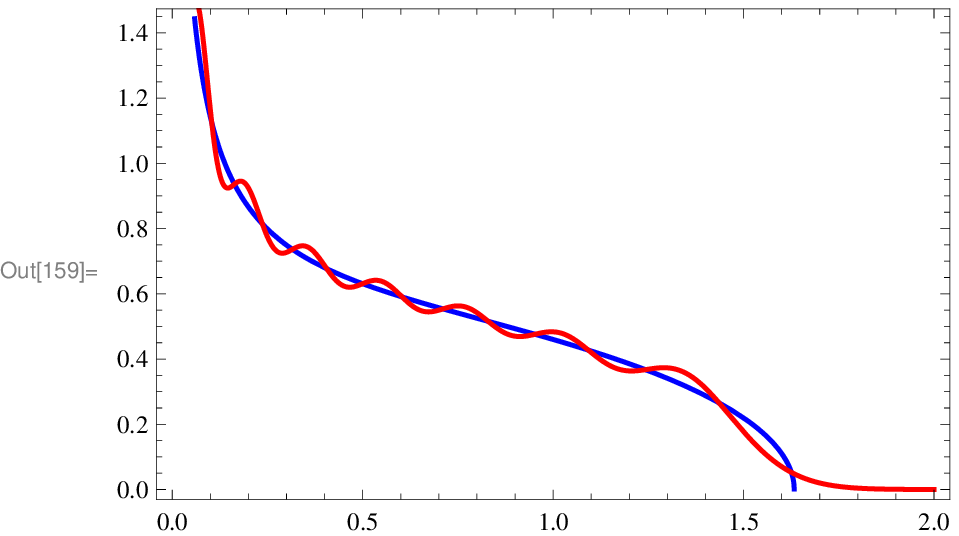}}

\end{figure}

\begin{flushleft}
\textcolor{blue}{\rule{1cm}{1pt}} Exact density $f_0$ of positive eigenvalues .\\
\textcolor{red}{\rule{1cm}{1pt}} Approximative density $f_7$ of positive eigenvalues.\\
\end{flushleft}

{\bf Second case:} $n=5, \mu_5=\frac52$, hence $\alpha=\frac12$.
$$\displaystyle f_{\frac12}(x)=\frac1{\pi}\sqrt{(1.9-x)(x-0.1)}\big(1+\frac1{2\sqrt{0.19}x}\big), {\rm with\;\; support}\;\; [0.1,1.9].$$
\newpage
\begin{figure}[h]
\centering\scalebox{0.5}{\includegraphics[width=17cm, height=12cm]{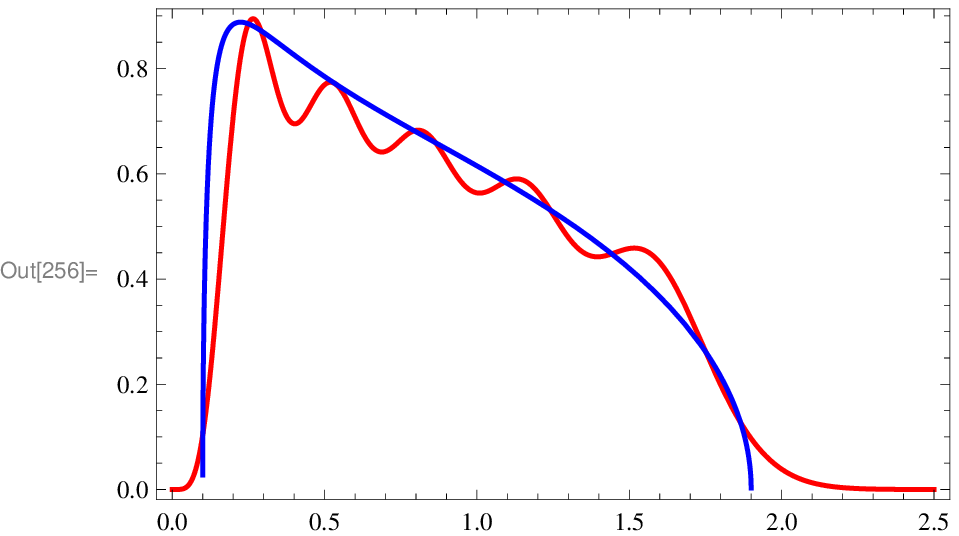}}

\end{figure}
\begin{flushleft}
\textcolor{blue}{\rule{1cm}{1pt}} Exact density $f_{\frac12}$ of positive eigenvalues .\\
\textcolor{red}{\rule{1cm}{1pt}} Approximative density $f_5$ of positive eigenvalues.\\
\end{flushleft}
\begin{center}Appendix\end{center}
\begin{lemma}\label{li} \
 For $z\in\Bbb C\setminus[a,b]$
$$\frac1{\pi}\int_a^b\frac1{z-t}\sqrt{\frac{t-a}{b-t}}dt=\sqrt{\frac{z-a}{z-b}}-1.$$
$$\frac1{\pi}\int_a^b\frac1{z-t}\sqrt{\frac{b-t}{t-a}}dt=-\sqrt{\frac{z-b}{z-a}}+1.$$
$$\frac1\pi\int_a^b\frac1{z-x}\sqrt{(b-x)(x-a)}dx=z-\frac{a+b}{2}-\sqrt{(z-a)(z-b)}.$$
\end{lemma}
{\it Consequence.}
$$\frac1\pi\int_a^b\sqrt{(b-x)(x-a)}dx=\frac{(a-b)^2}{8}.$$
$$\frac1{\pi}\int_a^b\sqrt{\frac{t-a}{b-t}}=\frac1{\pi}\int_a^b\sqrt{\frac{b-t}{t-a}}=\frac{b-a}{2}.$$
{\bf Proof.} The Cauchy-Stieljes transform
$$G(z)=\int_a^b\frac1{z-t}\sqrt{\frac{t-a}{b-t}}dt,$$ satisfies
$$[G]=-2i\pi\sqrt{\frac{t-a}{b-t}}\chi(t),$$
where for $x\in[a,b]$,\;\;$[G]=\displaystyle\lim_{\varepsilon\to 0}G(x+i\varepsilon)-G(x-i\varepsilon)$ is the bounded values distribution, and $\chi$ is the indicator function of $[a,b]$. Moreover
$$[\sqrt{\frac{z-a}{z-b}}]=-2i\sqrt{\frac{t-a}{b-t}}\chi(t).$$
The function $$\frac1{\pi}G(z)-\sqrt{\frac{z-a}{z-b}},$$
is holomorphic in $\Bbb C$ with limit $-1$ at infinity. By Liouville theorem one gets
$$\frac1{\pi}G(z)-\sqrt{\frac{z-a}{z-b}}=-1.$$
The proof of the second and third relation can be performed by the same lines with few modification.

\begin{lemma} For $a<b$
$$\int_a^bQ'(t)\frac{dt}{\sqrt{(t-a)(b-t)}}=\int_0^1Q'\big((1-t)b+ta\big)\frac{dt}{\sqrt {t(1-t)}}.$$
$$\int_a^bQ'(t)\sqrt{\frac {t-a}{b-t}}dt=(b-a)\int_0^1Q'\big((1-t)b+ta\big)\sqrt{\frac {1-t}{t}}dt.$$
\end{lemma}
{\bf Proof.} We use the change of variable $t=(1-u)b+au$.

For $a\in\Bbb R$, $b>a$, let $$\varphi(a,b)=\int_0^1Q'((1-t)b+ta)\frac{dt}{\sqrt{t(1-t)}},$$
 $$\psi(a,b)=(b-a)\int_0^1Q'((1-t)b+ta)\sqrt\frac {1-t}tdt-2\pi.$$
And $$h(a,b)=\frac{\psi(a,b)}{b-a}.$$
$Q$ a convex polynomial with degree $2m$, $m\geq 1$, and strictly positive leading coefficient.
\begin{proposition}\label{pp1} \
\begin{description}

\item[(1)] For all $a\in\Bbb R$ there is a unique $b(a)>a$, such that $\psi\big(a,b(a)\big)=0$.
\item[(2)] Let $\kappa=\inf\big\{a\in\Bbb R\mid \varphi(a,b(a))\geq 0\big\}$ which is an element of $[-\infty,+\infty[$. The map $a\mapsto b(a)$
define a increasing function on $]\kappa,+\infty[$,
and for all $a> \kappa$, $\varphi\big(a,b(a)\big)\geq 0$.
\item[(3)] There is a unique $a_0\in\Bbb R$, and a unique $b_0>a_0$, such that $\varphi(a_0,b_0)=0\;and\;\psi(a_0,b_0)=0,$ $(b_0=b(a_0))$.

    \end{description}
\end{proposition}
{\bf Proof.---}

(1)
For $a\in\Bbb R$, by the derivative theorem under the integral sign one gets
$$\frac{\partial h(a,b)}{\partial b}=\int_0^1Q''((1-t)b+ta)(1-t)\sqrt\frac {1-t}tdt+\frac{2\pi}{(b-a)^2}.$$
Since the polynomial $Q$ is convex, hence the function $b\mapsto h(a,b)$ increases strictly. Moreover, using Fatou's lemma and a simple computations give
$\displaystyle\lim_{b\to +\infty}h(a,b)=+\infty.$
From the dominate convergence theorem we obtain $\displaystyle\lim_{b\to a^+}h(a,b)=+\infty.$ Hence, the
 Rolle's theorem and the monotony of $h$ ensure the existence and uniqueness of $b(a)>a$ such that
$h(a,b(a))=0$ which means that $\psi(a,b(a))=0$.

(2) It is enough to prove that $\big\{a\in\Bbb R\mid \varphi(a,b(a))\geq 0\big\}\neq\emptyset.$
If for all $a\in\Bbb R$, $\varphi(a,b(a))<0$, then for all $a\in\Bbb R$,
 $$\int_0^1Q'((1-t)b(a)+ta)\frac{dt}{\sqrt{t(1-t)}}<0,$$
 since $b(a)>a$, hence
  $$\pi Q'(a)\leq \int_0^1Q'((1-t)b(a)+ta)\frac{dt}{\sqrt{t(1-t)}}<0.$$
  As $a\to +\infty$ one gets a contradiction. Thus $$\kappa=\inf\big\{a\in\Bbb R\mid \varphi(a,b(a))\geq 0\big\}\in[-\infty,+\infty[.$$
The map $a\mapsto b(a)$ is correctly defined from the unicity of the solution.

Furthermore for all $a> \kappa$, $$\varphi(a,b(a))\geq0.$$ follows from the definition of $\kappa$.

{\it Growth of $a\mapsto b(a)$.}
In the first hand by derivative theorem under the integral sign we have
$$\frac{\partial\psi(a,b)}{\partial a}=-\int_0^1Q'((1-t)b+ta)\sqrt{\frac {1-t}t}dt+(b-a)\int_0^1Q''((1-t)b+ta)\sqrt{t {(1-t)}}dt.$$
Using a integration by part in the second member, we obtain
$$\frac{\partial\psi(a,b)}{\partial a}=-\frac12\varphi(a,b).$$
Moreover $$\frac{\partial\varphi(a,b)}{\partial a}=\int_0^1Q''((1-t)b+ta)\sqrt{\frac t{1-t}}dt,$$
which is strictly positive by convexity of $Q$. Hence the function $\displaystyle a\mapsto\frac{\partial\psi(a,b)}{\partial a}$
  decreases strictly on $]-\infty,b[$ for all $b$. Then for all $(a_1,a)$ with  $\kappa<a_1<a$ and $a_1,a\in]-\infty,b(a_1)[$, one gets
$$\frac{\partial\psi(a,b(a_1))}{\partial a}<\frac{\partial\psi(a_1,b(a_1))}{\partial a}=-\frac12\varphi(a_1,b(a_1))\leq 0.$$
Hence the function $\displaystyle a\mapsto\psi(a,b(a_1))$ decreases strictly on $]-\infty,b(a_1)[$ .

Let $a_1<a_2$ and assume $b(a_2)<b(a_1)$, hence $a_1<a_2<b(a_2)< b(a_1)$. Moreover, we saw that for $b>a$,
the function $\displaystyle b\mapsto h(a,b)=\frac{\psi(a,b)}{b-a}$ is strictly increasing. Then
$$h(a_2,b(a_2))< h(a_2,b(a_1))=\frac{\psi(a_2,b(a_1))}{b(a_1)-a_2}\leq\frac{\psi(a_1,b(a_1))}{b(a_1)-a_2}=0,$$
thus $\psi(a_2,b(a_2)<0$, this give a contradiction. Which means for all $a_1< a_2$ then $b(a_1)< b(a_2)$ (here we used $b(a_1)\neq b(a_2)$ for $a_1\neq a_2$).

(3)
Let $$E=\Big\{a\in\Bbb R\mid \varphi(a,b(a))\geq 0\;and\;\psi(a,b(a))=0\Big\}.$$
Where $b(a)$ is the unique solution in $]a,+\infty[$ of the equation $\psi(a,b)=0$.
We want to show that $E$ admit a minimum. Indeed from the property (2), $E\neq\emptyset$.

Assume $E$ is unbounded below.

 From the monotony of $a\mapsto b(a)$, we have for all $a\leq 0$, $b(a)\leq b(0)$. Hence

 Then $\displaystyle\lim_{a\to-\infty}Q'((1-t)b(a)+t a)=-\infty$ for all $t\in[0,1]$, where we use the positivity of the leading coefficient of the polynomial $Q$ and the degree of $Q'$ is $2m-1$.

 Letting $a\to -\infty$ in the equation below
 $$-\varphi(a,b(a))=\int_0^1-Q'((1-t)b(a)+ta)\frac1{\sqrt{t(1-t)}}dt\leq 0,$$
 and using Fatou's lemma we get a contradiction. which prove that $E$ is bounded below. The closeness of $E$ is an immediate consequence of the continuity of the functions $(u,v)\mapsto\varphi(u,v)$ and $(u,v)\mapsto\psi(u,v)$ and the increasing of the function $a\mapsto b(a)$.

 We denote in the sequel $a_0=\min E$. From the previous we saw, there exist a unique $b_0:=b(a_0)>a_0$ such that $\varphi(a_0,b_0)\geq 0$ and $\psi(a_0,b_0)=0$.

For $n\in\Bbb N$, let $a_n=a_0-\frac 1n$, Since $a_n\notin E$, hence for all $b\geq a_n$,
$\varphi(a_n, b)<0$ or  $\psi(a_n, b)\neq 0$. Moreover for all $n\in\Bbb N$,
from property (1) of the proposition there exist a unique $b_n>a_n$,
such that $\psi(a_n, b_n)=0$, thus $\varphi(a_n,b_n)<0$.
Since the sequence $a_n$ converges to $a_0$ and $b_n>a_n$, then $b_n$ is bounded from below.
If $b_n$ is unbounded from above, then $\displaystyle\lim_{n\to +\infty}b_n=+\infty$.
 By continuity of the function $\varphi$ and Fatou's lemma we have $\displaystyle\lim_{n\to +\infty}\varphi(a_n,b_n)=+\infty\leq 0$ this give a contradiction. Hence the sequence $b_n$ is bounded.

 Let $y_n$ be a convergent subsequence of $b_n$ with limit $y\geq a_0$. We have
$$\psi(a_n, y_n)=0\;and\;\varphi(a_n,y_n)\leq 0.$$
By continuity of $\varphi(u,v)$ and $\psi(u,v)$, as $n$ goes to $+\infty$ we obtain
$$\psi(a_0, y)=0\;and\;\varphi(a_0,y)\leq 0.$$
Since the equation $\psi(a_0,y)=0$ has a unique solution $b_0$ in $]a_0,+\infty[$ and $\psi(a_0,a_0)=-2\pi$ hence $y=b_0$.
Thus, $$\varphi(a_0,b_0)\leq 0.$$
Together with the reverse inequality one gets the desired result.

\begin{center}{\bf Acknowledgments.}
\end{center}
My sincere thanks go to Jacques Faraut for his comments on this manuscript and his important remarks.

Address:  College of Applied Sciences
   Umm Al-Qura University
  P.O Box  (715), Makkah,
  Saudi Arabia.\\
  Facult\'e des Scinces de Tunis,  Campus Universitaire El-Manar, 2092 El Manar Tunis.\\
E-mail: bouali25@laposte.net \& mabouali@uqu.edu.sa

\end{document}